\documentclass[onefignum,onetabnum]{siamart190516}

\usepackage{amsfonts}
\usepackage{graphicx}
\usepackage{epstopdf}
\usepackage{booktabs}
\usepackage{hyperref}
\usepackage{algorithmic}
\ifpdf
  \DeclareGraphicsExtensions{.eps,.pdf,.png,.jpg}
\else
  \DeclareGraphicsExtensions{.eps}
\fi
\usepackage{pgfplots}
\pgfplotsset{compat=1.14}
\usepgfplotslibrary{colorbrewer}

\pgfplotsset{every axis/.append style={
            		width=2.5in,
            		grid=both,
                    label style={font=\footnotesize},
                    tick label style={font=\footnotesize},
            		cycle list/Set1-5,
            		cycle multiindex* list={
            		    Set1-5
            		    \nextlist
                        marks
                        \nextlist
                        [2 of]linestyles
                        \nextlist
            		    thick
            		    \nextlist
                    },
            }
}

\pgfplotsset{
    cycle list/.define={marks}{
	    every mark/.append style={solid,fill opacity=0.25,fill=\pgfkeysvalueof{/pgfplots/mark list fill}},mark=*\\
	    every mark/.append style={solid,fill opacity=0.25,fill=\pgfkeysvalueof{/pgfplots/mark list fill}},mark=square*\\
	    every mark/.append style={solid,fill opacity=0.25,fill=\pgfkeysvalueof{/pgfplots/mark list fill}},mark=triangle*\\
	    every mark/.append style={solid,fill opacity=0.25,fill=\pgfkeysvalueof{/pgfplots/mark list fill}},mark=diamond*\\
    },
}

\newcommand{\norm}[1]{\lVert #1 \rVert}

\headers{An Iterative Method For Contour Based Nonlinear Eigensolvers}{J. Brenneck, E. Polizzi}

\title{An Iterative Method For Contour Based Nonlinear Eigensolvers\thanks{Submitted to the editors 6/30/20.
\funding{National Science Foundation grants SI2-SSE-1739423 and AF-1813480.}
}}

\author{Julien Brenneck\thanks{Department of Mathematics and Statistics, University of Massachusetts, Amherst
	(\email{jbrenneck@math.umass.edu}, \email{polizzi@ecs.umass.edu}).}
  \and Eric Polizzi\footnotemark[2]
}

\usepackage{amsopn}

\ifpdf
\hypersetup{
  pdftitle={An Iterative Method for Contour Based Nonlinear Eigensolvers},
  pdfauthor={J. Brenneck, E. Polizzi}
}
\fi

\begin{document}

\maketitle

\begin{abstract}
Contour integration techniques have become a popular choice for solving the linear and non-linear eigenvalue problems.
  They principally include the Sakurai-Sugiura methods, the Beyn's algorithm, the FEAST/NLFEAST algorithms and  other rational filtering techniques.
  While these methods can result in effective 'black-box' approach for solving linear eigenvalue problems,
  they still present several shortcomings for addressing nonlinear eigenvalue problems which are both mathematically and practically far more challenging.
  In this paper, we introduce a new hybrid algorithm that advantageously combines the iterative nature of NLFEAST with the effectiveness of Beyn's approach to deal with general non-linearity.
  In doing so, this NLFEAST-Beyn hybrid algorithm can overcome current limitations of both algorithms taken separately.
  After presenting its derivation from both a Beyn's and NLFEAST's perspective, several numerical examples are discussed
  to demonstrate the efficiency of the new technique.
\end{abstract}

\begin{keywords}
  nonlinear eigenvalue problem, contour integration, FEAST, Beyn's algorithm, residual inverse iteration
\end{keywords}

\begin{AMS}
  15A18, 65H17
\end{AMS}

\section{Introduction}
Eigenvalue problems in which the coefficient matrices depend nonlinearly on the eigenvalues arise in a variety of applications in science and engineering including dynamic analysis of structures, vibrations of fluid-solid structures, and computational nanoelectronics, to name just a few. 
A large collection of examples has been compiled in NLEVP~\cite{betcke2013}, as well as an associated MATLAB package. 
The resulting nonlinear eigenvalue problem (NEP) takes the form:
\begin{equation}
	T(\lambda) x = 0, \quad
	x \in \mathbb{C}^n, \quad
	\lambda \in \Omega
        \label{eq:nl}
\end{equation}
where the nonlinear function \( T \in H(\Omega, \mathbb{C}^{n\times n}) \) is holomorphic on some domain \( \Omega \subset \mathbb{C} \).
This is a numerically challenging problem, which includes the linear eigenvalue problem as a special case, letting \( T(z) = A - z I \), as well as the rootfinding problem, for example letting \( T(z) = \cos(z) \).
There are practical difficulties associated with the nonlinear problem not seen in the linear case, in particular there can be infinitely many eigenvalues, as well as linearly dependent eigenvectors~\cite{voss2014}.
We refer to the review article by G\"uttel and Tisseur~\cite{guttel2017} for an overview of the nonlinear eigenvalue problem, which includes a section on contour based methods.

In most applications, one is interested in finding all the eigenvalue-eigenvector pairs in a given region of the complex plane.
Contour based methods are well suited for this problem, since they aim at finding all eigenpairs \( (x_i, \lambda_i) \) solving (\ref{eq:nl}) associated with the eigenvalues inside some contour \( \Gamma \subset \Omega \). Solving within a well chosen contour reduces the problem to a more tractable one. 
Critically, the size of the resulting reduced problem is largely independent of the size of the original problem.
The approximation of the contour by numerical quadrature, as well as the ability to arbitrarily split regions of the spectrum into independent problems, gives these methods an inherently parallel structure, a necessity for the increasingly distributed nature of computational hardware.
The contour based approach is ultimately based on the Cauchy integral formula, which in matrix form~\cite{higham2008} is given as
\begin{equation}\label{eq:cauchy}
	\frac{1}{2\pi i} \int_{\Gamma} f(z)T^{-1}(z) \, dz = V f(J) W^H.
\end{equation}
Here we assume \( T(z) \) has a finite number of  eigenvalues in the contour \(\Gamma\), and \(f(z)\) is holomorphic, then \( J \) is a block Jordan matrix with left and right generalized eigenvectors \( W \) and \( V \), respectively.

The previous contour based NEP solvers considered in this paper are Beyn's method~\cite{beyn2012,vanbarel2016}, the SS methods~\cite{asakura2009,yokota2013}, and NLFEAST~\cite{gavin2018,gavin2018a}. 
Other contour based algorithms are also possible using rational filtering techniques \cite{saad2019,vanbarel2016b}.
The approach of Beyn's method, as well as the method we present, probe the 
Jordan decomposition in (\ref{eq:cauchy}) for spectral information.



This paper introduces a novel algorithm based largely on Beyn's method and the NLFEAST algorithm.
While Beyn's algorithm offers simplicity and effectiveness in dealing with general non-linear problem, it is not iterative in nature and requires expensive numerical contour integration with a large number of quadrature nodes (each node corresponds to a linear system solve). 
On the other hand, NLFEAST benefits from a highly efficient iterative procedure that can systematically reach convergence without using an excessive number of quadrature nodes per iteration (with the possibility to also use low-accuracy and precision for the system solves).
NLFEAST, however, leads to a reduced eigenvalue problem which is still nonlinear in nature and cannot be easily solved beyond the polynomial form.
By combining the Beyn and NLFEAST approaches, we propose a new hybrid algorithm that aims at addressing their own separate shortcomings.
The paper is organized as follows. An overview of the previous algorithms is given in \cref{sec:beyn},
our new algorithm is presented in \cref{sec:alg} and \cref{sec:moments}, and numerical
results are discussed in \cref{sec:experiments}.

\section{Background: the Beyn and FEAST Algorithms}
\label{sec:beyn}

We first give some definitions necessary to state the algorithms, as well as the Keldysh theorem. 

The function \( T(z) \) is \textit{regular} if \( \det T(z) \) does not vanish identically on \( \Omega \), or equivalently the resolvent set \( \rho(T) = \Omega \setminus \sigma(T) \) is non-empty, where \( \sigma(T) \) is the spectrum.
We assume throughout this paper that \( T(z) \) is regular on a non-empty domain \( \Omega \subset \mathbb{C} \).

\begin{definition}
	Let \( T \in H(\Omega, \mathbb{C}^{n\times n}) \) and \( \lambda \in \Omega \). 
	\begin{enumerate}
		\item	
			A vector valued function \( v \in H(\Omega, \mathbb{C}^n) \) is called a \textit{root function of \( T \) at \( \lambda \)} if
			\begin{equation*}
				v(\lambda) \neq 0, 
				\quad
				T(\lambda)v(\lambda) = 0.
			\end{equation*}
			The multiplicity of the root \( z = \lambda \) of \( T(z)v(z) \) is denoted \( s(v) \).
		\item
			A tuple \( (v_0, \cdots v_{m-1}) \in {(\mathbb{C}^n)}^m \) with \( m \ge 1 \) and \( v_0 \neq 0 \) is called a \textit{Jordan chain for \( T \) at \( \lambda \)} if \( v(z) = \sum_{k=0}^{m-1} {(z-\lambda)}^k v_k \) is a root function for \( T \) at \( \lambda \) and \( s(v) \ge m \).
		\item
			For a given \( v_0 \in N(T(\lambda)) \), and \( v_0 \neq 0 \), the number
			\begin{equation*}
				r(v_0) = \max \{ s(v) : v \text{ is a root function for } T \text{ at } \lambda \text{ with } v(\lambda) = v_0 \}
			\end{equation*}
			is finite and called the \textit{rank of \( v_0 \)}.
		\item
			A system of vectors in \( \mathbb{C}^n \),
			\begin{equation*}
				V = (v_k^j : 0 \le k \le m_j -1, \; 1 \le j \le d),
			\end{equation*}
			is a \textit{complete system of Jordan chains for \( T \) at \( \lambda \)} if
			\begin{enumerate}
				\item
					\( d = \dim N(T(\lambda)) \) and \( \{v_0^1, v_0^2, \ldots, v_0^d \} \) is a basis for \( N(T(\lambda)) \).
				\item
					The tuple \( (v_0^1, \ldots, v_{m_j-1}^j) \) is a Jordan chain for \( T \) at \( \lambda \) for \( j=1,\ldots,d \).
				\item
					\( m_j = \max \{r(v_0) : v_0 \in N(T(\lambda)) \setminus \text{span} \{v_0^\nu : 1 \le \nu < j\} \} \) for \( j = 1,\ldots,d \).
			\end{enumerate}
	\end{enumerate}
\end{definition}
	One can show that a complete system of Jordan chains always exists, and that the numbers \( m_j \) satisfy \( m_1 \le m_2 \le \cdots \le m_d \), they are called the \textit{partial multiplicities} of \( \lambda \).
	We can now state Keldysh's theorem, and we do so in a concise matrix notation~\cite{guttel2017}.

	\begin{theorem}[Keldysh]
		Let \( T \in H(\Omega, \mathbb{C}^{n\times n}) \) and \( \lambda_1, \ldots, \lambda_s \) be the distinct eigenvalues of \( T \) in \( \Omega \) of partial multiplicities \( m_{i,1} \ge \cdots \ge m_{i,d_i} \) and define
		\begin{equation*}
			\overline{m} = \sum_{i=1}^s \sum{j=1}^{d_i} m_{ij}
		\end{equation*}
		then there are \( n \times \overline{m} \) matrices \( V \) and \( W \) whose columns are generalized eigenvectors, and an \( \overline{m} \times \overline{m} \) Jordan matrix \( J \) with eigenvalues \( \lambda_i \) of partial multiplicities \( m_{ij} \), such that
		\begin{equation*}
			{T(z)}^{-1} = V {(zI - J)}^{-1} W^H + R(z)
		\end{equation*}
		for some \( R \in H(\Omega, \mathbb{C}^{n\times n}) \).
	\end{theorem} 
	There are two important special cases~\cite{guttel2017}.
	When all eigenvalues \( \lambda_i \) are semisimple, the matrix \( J \) is diagonal.
	In the case when all eigenvalues \( \lambda_i \) are simple, the matrix \( J \) is diagonal and \( V, W \) consist of right and left eigenvectors satisfying \( w_i^H T'(\lambda_i) v_i = 1 \).

	We now consider a simple closed contour \( \Gamma \subset \Omega \), and denote by \( n(\Gamma) \) the number of eigenvalues of \( T \) inside the interior of \(  \Gamma \). 
	The result (\ref{eq:nl}) is then obtained from the residue theorem.
    We note that the function \( f(J) \) represents the standard application of a function to a Jordan matrix.
    When \( J \) is diagonal, this reduces to mapping \( f \) onto the diagonal elements.

    The various contour based algorithms, including Beyn's method, NLFEAST, and the SS methods, use the moments of the Cauchy integral of the resolvent applied to a probing matrix (or initial subspace) \( X \in \mathbb{C}^{n\times m} \) with subspace size \(m\geq n(\Gamma)\), typically this is taken to be random.
    The moments are then defined by
\begin{equation}\label{eq:moments}
	A_k = \frac{1}{2 \pi i} \int_{\Gamma} z^k T^{-1}(z) X \, dz.
\end{equation}
Beyn's method typically uses the zeroth and first moments (\( A_0, A_1 \)), and can be extended to use higher moments for problems with many eigenvalues in the contour or where there are linearly dependent eigenvectors.
The SS methods uses as many moments as are necessary to resolve the eigenvalues within the contour, whereas NLFEAST uses only the zeroth moment as a projector and converges through residual inverse iteration.
The Beyn and SS methods are not iterative, and convergence is dependant on the accuracy of the quadrature, though the SS-type methods can be viewed as iteratively increasing the subspace by adding moments until the solution is sufficiently accurate~\cite{saad2019}. 

\subsection{Beyn's Method}

There are two algorithms presented by Beyn~\cite{beyn2012}, we consider here the first algorithm for simple eigenvalues where the eigenvectors are linearly independant. Beyn's method applies the theorem of Keldysh, giving a linearization of the problem in terms of the moments \( A_0 \) and \( A_1 \) that can be solved using traditional linear solvers.
In the case of simple eigenvalues this gives
\begin{equation}
	A_0 = V W^H X,
	\quad\quad\quad
	A_1 = V \Lambda W^H X
\end{equation}
leading to a rectangular equation of the form \(A_1x = A_0x \lambda\).
Using a singular value decomposition of \( A_0 = V_0 \Sigma_0 W_0^H \) one can derive a computable matrix similar to \( \Lambda \).
\begin{equation}\label{eq:beyn}
	\Lambda \sim (V_0^H V) \Lambda {(V_0^H V)}^{-1} = V_0^H A_1 W_0 \Sigma_0^{-1} 
\end{equation}
Letting \( B = V_0^H A_1 W_0 \Sigma_0^{-1} \) we can solve the linear standard eigenvalue problem \( BY = Y\Lambda \) and obtain solutions to the original problem by taking \( X = V_0 Y \).
The full algorithm, with the resizing step based on singular values removed for clarity, is presented in algorithm~\ref{alg:beyn}. 

\begin{algorithm}
\caption{Beyn's Method}
\label{alg:beyn}
\begin{algorithmic}
	\REQUIRE{Initial (random) subspace \( X \in \mathbb{C}^{n \times m}\)}
	\REQUIRE{Contour \( \Gamma \) and \( N \) quadrature nodes and weights \( (z_j, \omega_j) \)}
	\STATE{\(A_0 = \sum_{j=1}^{N} \omega_j {T(z_j)}^{-1} X \)}
	\STATE{\(A_1 = \sum_{j=1}^{N} \omega_j z_j {T(z_j)}^{-1} X \)}
	\STATE{Compute the Singular Value Decomposition \(V_0 \Sigma_0 W_0^H \leftarrow A_0\)}
	\STATE{\(B = V_0^H A_1 V_0 \Sigma_0^{-1}\)}
	\STATE{Solve \( BY = Y\Lambda \)}
	\STATE{ \( X \leftarrow V_0 Y \)}
	\RETURN \( \Lambda, X \)
\end{algorithmic}
\end{algorithm}

There are practical limitations to Beyn's method, in particular due to the numerical quadrature. 
When using Beyn's method to solve a problem, one does not know how many quadrature nodes will be necessary for sufficiently accurate solutions.
Without an adaptive quadrature scheme, inaccurate solutions means Beyn's method must be run again with a more accurate integration.
Each quadrature node corresponds to a factorization and linear system solve in the algorithm, all of which can be done in parallel.
In practice, however, high accuracy may require a large number of quadrature nodes (typically greater than $64$) which is numerically expensive.
Furthermore, a high number of contour nodes may increase the chances to generate highly ill-conditioned linear systems (e.g. if a node ends up too close to an eigenvalue), causing numerical issues, as probing a Jordan decomposition is highly sensitive to perturbations~\cite{guttel2017}.
One of the aims of this paper is to present a method of iterative refinement for Beyn's method, giving convergence with a small number of quadrature nodes.

\subsection{NLFEAST}

The original FEAST algorithm \cite{polizzi2009,tang2014} can be interpreted as a generalization of the shift-and-invert iteration for solving the linear eigenvalue problem using multiple fixed shifts along a contour integration.
The residual inverse iteration (RII)  proposed by Neumaier~\cite{neumaier1985}  is a modification of the shift-and-invert iteration which allows it to be applied to nonlinear problems.
The key insight of the NLFEAST algorithm~\cite{gavin2018,gavin2018a} is to generalize RII using multiple fixed shifts, and similarly to FEAST using the quadrature nodes of the contour integration as shifts.
We define the block form of the residual function for \( X \in \mathbb{C}^{n \times m} \) as
\begin{equation}
	T(X, \Lambda) = \big[ T(\lambda_1)x_1, T(\lambda_2)x_2, \; \cdots \;, T(\lambda_m)x_m \big]
\end{equation}
This gives the following way of computing the zeroth moment.
\begin{equation}\label{eq:nlfeast_contour}
	Q_0 = \frac{1}{2 \pi i} \int_{\Gamma} \big( X - T^{-1}(z)T(X, \Lambda)\big) {(zI - \Lambda)}^{-1} \, dz
\end{equation}
For \( T \) a standard linear problem we can show that \( Q_0 = A_0 \), but this does not hold in general~\cite{neumaier1985}, thus the moments used in the computation are fundamentally different.
The approach of NLFEAST is to use \( Q_0 \) as a projector, which leads to the following reduced nonlinear problem:
\begin{equation}\label{eq:red}
	Q_0^H T(\lambda) Q_0 y = 0	
\end{equation}
from which one can recover eigenpairs \( (\lambda, Q_0 y) \).
This resulting nonlinear eigenvalue problem of reduced dimension must in turn be solved using any suitable method.
A modified version of NLFEAST, reordered with an explicit first iteration for clarity, is presented in algorithm~\ref{alg:nlfeast}.
While RII itself is well studied~\cite{jarlebring2011}, the multiple shift contour approach of NLFEAST and convergence properties are not yet fully understood.
Analogously to the linear FEAST algorithm, numerical experiments show that NLFEAST can successfully converge using a relatively small number of contour nodes, and that the convergence rate can be systematically improved using larger search subspaces or by adding additional contour nodes \cite{gavin2018}.
One of the most attractive feature of RII is that it allows the linear systems to be solved with low accuracy
(e.g. using single precision arithmetic and/or inexact iterative solves) while preserving, for the most part, the FEAST convergence rate~\cite{gavin2018a,golub2000a}.
RII has then become the new {\em de facto} standard to all FEAST algorithms implemented in the FEAST numerical library \cite{polizzi2020}.


\begin{algorithm}
\caption{NLFEAST}
\label{alg:nlfeast}
\begin{algorithmic}
	\REQUIRE{Initial (random) subspace \( X \in \mathbb{C}^{n \times m_0}\)}
	\REQUIRE{Contour \( \Gamma \) and \( N \) quadrature nodes and weights \( (z_j, \omega_j) \)}
	\REQUIRE{Stopping tolerence \( \epsilon \)}
	\STATE{ Orthogonalize \(Q = \sum_{j=1}^{N} \omega_j  T^{-1}(z_j) X \) }
	\STATE{Solve the nonlinear equation \( Q^H T(\lambda) Q y = 0 \) for approximate eigenpairs \( (\lambda, Qy ) \)}
	\STATE{\( \Lambda, X \leftarrow \text{diag}(\lambda_1, \ldots, \lambda_{m_0}), \;  [Qy_1, \ldots, Qy_{m_0}] \)}
	\WHILE{not converged}
	\STATE{Orthogonalize \(Q \leftarrow \sum_{j=1}^{N} \omega_j \big(X - T^{-1}(z_j) T(X, \Lambda)\big) {(z_j I - \Lambda)}^{-1} \) }
	\STATE{Solve the nonlinear equation \( Q^H T(\lambda) Q y = 0 \) for approximate eigenpairs \( (\lambda, Qy ) \)}
	\STATE{\( \Lambda, X \leftarrow \text{diag}(\lambda_1, \ldots, \lambda_{m_0}), \;  [Qy_1, \ldots, Qy_{m_0}] \)}
	\ENDWHILE
\end{algorithmic}
\end{algorithm}

The new FEAST v4.0 package implements a simple linearization of the reduced system~(\ref{eq:red})
(i.e.~by forming the companion problem) for solving the polynomial eigenvalue problem~\cite{polizzi2020}. For addressing the general non-linear problem,
we have succesfully been testing with the Beyn's algorithm for solving~(\ref{eq:red})~\cite{gavin2018a}. 
In practice, however, making use of two levels of contour integration
(at the level of both the original and the reduced systems)
leads to an expensive numerical procedure, and may appear somehow redundant.
The goal of this paper is to address this issue using a new more effective NLFEAST-Beyn hybrid approach.

\section{An NLFEAST-Beyn Hybrid Algorithm}
\label{sec:alg}

We develop a new algorithm by applying the residual inverse iteration to the complex moments of~(\ref{eq:moments}), generalizing the contour integration in ~(\ref{eq:nlfeast_contour}) of NLFEAST.
\begin{equation}
	Q_k = \frac{1}{2 \pi i} \oint_{\Gamma} z^k \Big(X - T^{-1}(z)T(X, \Lambda) \Big){(zI - \Lambda)}^{-1} \, dz
\end{equation}
We then apply the linearization technique of Beyn's method to these moments.
In NLFEAST this linearization would be applied to the reduced problem after a projection using Rayleigh-Ritz, whereas here we apply it directly.
In the linear case we have~\cite{neumaier1985} that \( A_k = Q_k \), though this does not hold for \( T(\lambda) \) nonlinear.
The new method is given in algorithm~\ref{alg:hybrid}, presented with a QR approach instead of SVD to highlight similarities with NLFEAST.
In particular, the existing NLFEAST code is easily adapted into this form.

\begin{algorithm}
\caption{NLFEAST-Beyn Hybrid Algorithm}
\label{alg:hybrid}
\begin{algorithmic}
	\REQUIRE{Initial (random) subspace \( X \in \mathbb{C}^{n \times m}\)}
	\REQUIRE{Contour \( \Gamma \) and \( N \) quadrature nodes and weights \( (z_j, \omega_j) \)}
	\REQUIRE{Stopping tolerence \( \epsilon \)}
	\STATE{\(Q_0 = \sum_{j=1}^{N} \omega_j {T(z_j)}^{-1} X \)}
	\STATE{\(Q_1 = \sum_{j=1}^{N} \omega_j z_j {T(z_j)}^{-1} X \)}
	\STATE{Compute the QR Decomposition \( qr \leftarrow Q_0\)}
	\STATE{\( B = q^H Q_1 r^{-1} \)}
	\STATE{Solve \( BY = Y \Lambda \)}
	\STATE{ \( X \leftarrow qY \)}
	\WHILE{not converged}
\STATE{\(Q_0 \leftarrow \sum_{j=1}^{N} \omega_j [X - {T(z_j)}^{-1}T(X, \Lambda)] {(z_j I - \Lambda)}^{-1} \)}
\STATE{\(Q_1 \leftarrow \sum_{j=1}^{N} \omega_j z_j [X  - {T(z_j)}^{-1} T(X, \Lambda)] {(z_j I - \Lambda)}^{-1} \)}
	\STATE{Compute the QR Decomposition \( qr \leftarrow Q_0\)}
	\STATE{\( B \leftarrow q^H Q_1 r^{-1} \)}
	\STATE{Solve \( BY = Y \Lambda \)}
	\STATE{ \( X \leftarrow qY \)}
	\ENDWHILE
	\RETURN \( X, \Lambda \)
\end{algorithmic}
\end{algorithm}

The NLFEAST and Beyn algorithms give practical methods for solving nonlinear eigenproblems within a contour with a similar numerical quadrature approach.
As previously discussed, NLFEAST must use some other technique internally for solving the reduced nonlinear system.
When Beyn's method is used within NLFEAST, the question arose of how the contour, and thus the solutions of the linear systems, could be shared between them, as this is the most expensive step in both algorithms.
From the other direction, the question was how Beyn's method could be iterated.
The following algorithm answers both of these.
From one perspective, it can be viewed as NLFEAST using the linearization (from Keldysh) of Beyn's method directly.
From the other, it can be viewed as Beyn's method using the residual inverse iteration (from Neumaier) approach of NLFEAST.

A theoretical benefit of this new algorithm is that it can be reduced to standard FEAST for linear problems, which is not true of NLFEAST, and suggests that from a theory perspective this may be a more natural extension of the FEAST type algorithms.
To see this equivalence we must use a more direct Rayleigh-Ritz approach, which is less useful numerically.
We consider a standard linear problem \( Ax = \lambda x \) with \( A \) unitarily diagonalizable, with decompositions \( Q_0 = VV^H X \) and \( Q_1 = V \Lambda V^H X \).
Then the Rayleigh-Ritz projection of linear FEAST solves
\begin{equation}\label{eq:feast_eig_prob}
	Q_0^H A Q_0 x = \lambda Q_0^H Q_0 x.
\end{equation}
Note that for the linearization of Beyn's method, we equivalently must solve the rectangular equation \( Q_1 x = \lambda Q_0 x \).
Applying a projection directly with \( Q_0^H \) gives
\begin{equation}
	Q_0^H Q_1 x = Q_0^H A Q_0 x = \lambda Q_0^H Q_0 x
\end{equation}
by noting that \( AQ_0 = Q_1 \), and we see this is equivalent to~(\ref{eq:feast_eig_prob}).
Thus the algorithms can be seen as making equivalent Rayleigh-Ritz projections in the linear hermitian case, though numerical issues result in both methods using different projections in practice.
This can be further extended to the generalized problem by instead projecting with \( Q_0^H B \), which is then equivalent to generalized FEAST using B-orthogonality. A two-contour strategy equivalence that computes left and right eigenvectors
for the linear non-hermitian problem~\cite{kestyn2016} is also possible. 


Beyn uses an SVD to give a computable similarity for \( \Lambda \), which has the benefits of somewhat improved numerical stability and information about the rank.
In particular, filtering the singular values enforces rank conditions, the importance of which is emphasized by Beyn~\cite{beyn2012}.
Here we present instead a similarity transform based on a QR decomposition, which in testing has resulted generally in only a small decrease of convergence speed.
The choice of decomposition should depend on the application, in particular rank information may be desireable.
Typically the linear system solves are the computational bottleneck, so the cost of the decomposition may be less important.
This is not always the case, for example if the factorizations of the linear systems can be stored in memory, which gives a significant improvement in runtime, the cost of the decomposition could dominate successive iterations.

\section{Higher Moments}
\label{sec:moments}
When there are linearly dependant eigenvectors it is necessary to use higher moments of the integral in~(\ref{eq:moments}).
We give a contructed example of such a problem in~\autoref{sec:experiments}, though such a situation should be expected when there are in general more eigenvectors than the dimension of the system.
Beyn's method employs higher moments also to handle when there are more eigenvalues in the contour than the dimension \( n \) of the system. 
The SS methods are based entirely on using higher order moments.
A general formulation encompassing both Beyn and the SS methods is presented here~\cite{guttel2017}.

In order to generalize both Beyn and SS-Hankel we must use a left probing matrix \( \widehat{X} \in \mathbb{C}^{n \times \ell} \) in addition to the right probing matrix \( X \in \mathbb{C}^{n \times m} \), giving
\begin{equation}\label{eq:moments_left}
	\widehat{A}_p = \widehat{X}^H A_p = \frac{1}{2 \pi i} \int_{\Gamma} z^p \widehat{X}^H {T(z)}^{-1} X \, dz.
\end{equation}
We choose \( K \in \mathbb{N} \) as the number of moments.
Then the moment matrices are used to form the \( K\ell \times K m \) block Hankel matrices
\begin{equation}
	H_0 = \begin{bmatrix}
		\widehat{A}_0 & \cdots & \widehat{A}_{K-1}  \\
		\vdots & & \vdots \\
		\widehat{A}_{K-1} & \cdots & \widehat{A}_{2K-2}  \\
	\end{bmatrix},
	\quad
	H_1 = \begin{bmatrix}
		\widehat{A}_1 & \cdots & \widehat{A}_{K}  \\
		\vdots & & \vdots \\
		\widehat{A}_{K} & \cdots & \widehat{A}_{2K-1}  \\
	\end{bmatrix}.
\end{equation}
Following from the Keldysh theorem we have \( \widehat{A}_p = \widehat{X} V\Lambda^p W^H X \), and writing
\begin{equation}
	V_{[K]} = \begin{bmatrix}
		\widehat{X}V  \\
		\vdots \\
		\widehat{X}V \Lambda^{K-1}
	\end{bmatrix}
	,\quad
	W_{[K]}^H = [W^H X, \cdots, \Lambda^{K-1} W^H X]
\end{equation}
we have representations for \( B_0, B_1 \) that give a similarity to \( \Lambda \).
\begin{equation}
	H_0 = V_{[K]} W_{[K]}^H,
	\quad
	H_1 = V_{[K]} \Lambda W_{[K]}^H
\end{equation}
Then we let \( H_0, H_1 \) replace \( A_0, A_1 \) in Beyn's method as formulated above, and taking the similarity transform in~\ref{eq:beyn} we can solve the NEP.
In particular we take the SVD \( H_0 = V_0\Sigma_0 W^H_0  \) and then diagonalize \( V_0^H H_1 W_0 \Sigma_0^{-1} \).
This gives a method for computing the eigenvalues, but recovering the eigenvectors depends on the form of \( \widehat{X} \).
For Beyn's method, we take \( \widehat{X} = I \in \mathbb{C}^{n\times n} \), and thus \( \widehat{A}_p = A_p \) and \( \ell = n \), allowing us to recover the eigenvectors as before.
Taking \( \ell = m  \) gives the block SS-Hankel method, which is solved using a different decomposition to recover the eigenvectors.

For this paper we consider only \( \widehat{X} = I \), as in Beyn's method.
The SS-Hankel approach is in some ways more attractive for large systems, as Beyn's method scales poorly in the amount of memory needed for large \( n, K \), exploring this is a potential direction for future work.
We present such a method in algorithm~\ref{alg:generalized}, which is stated generally, but implemented for numerical experiments using the SVD similarity approach of Beyn.

\begin{algorithm}
\caption{Hybrid Algorithm with Higher Moments}
\label{alg:generalized}
\begin{algorithmic}
	\REQUIRE{Initial (random) subspace \( X \in \mathbb{C}^{n \times m}\)}
	\REQUIRE{Contour \( \Gamma \) and \( N \) quadrature nodes and weights \( (z_j, \omega_j) \)}
	\REQUIRE{Stopping tolerance \( \epsilon \), number of moments \( K \)}
	\STATE{\(Q_k = \sum_{j=1}^{N} \omega_j z_j^k {T(z_j)}^{-1} X \)}
	\STATE{Compute \( H_0, H_1 \)}
	\STATE{Solve (e.g.~by SVD) for \( X, \Lambda \)}
	\WHILE{not converged}
	\STATE{\(Q_k \leftarrow \sum_{j=1}^{N} \omega_j z_j^k [X  - {T(z_j)}^{-1} T(X, \Lambda)] {(z_j I - \Lambda)}^{-1} \)}
	\STATE{Compute \( H_0, H_1 \)}
	\STATE{Solve (e.g.~by SVD) for \( X, \Lambda \)}
	\ENDWHILE
	\RETURN \( X, \Lambda \)
\end{algorithmic}
\end{algorithm}

In this way it is possible to apply our iterative method to the higher order moments, simply replacing \( A_k \) with analogous \( Q_k \) in the computation of \( B_0 \) and \( B_1 \).
The main difficulty in doing this is that it becomes necessary to deflate the subspace after every iteration, as we go from \( m \) vectors to \( K \times m \) vectors.
This inhibits such a method from computing \( m > n \) eigenvalues in a contour, a serious limitation for small highly nonlinear problems.
Selecting eigenpairs based on residual and distance from the center of the contour works well in practice. 
Whether the moment expanded subspace can be incorporated into the RII is a topic for future work.

\section{Experimental results}
\label{sec:experiments}

\begin{table}[htbp]
	\begin{center}
	\footnotesize
	\caption{Parameters of numerical experiments.}\label{tab:parameters}
	\begin{tabular}{lrrrrrrr}
		\toprule
		Problem & Type & Eigvals & \(K\) & Dim. \( (n) \) & \( m \) & Radius & Center \\
		\midrule
		\verb$butterfly$      & quartic     & 13 & 1 & 64   & 30 &   0.5 & \( 1+1i \) \\
		\verb$test_deficient$ & quadratic   & 4  & 2 & 15   & 4  &  0.25 & 0  \\
		\verb$gun$            & irrational  & 17 & 1 & 9956 & 32 & 30000 & 140000  \\
		\verb$hadeler$        & exponential & 12 & 1 & 200  & 15 &    10 & \(-30\)  \\
		\bottomrule
	\end{tabular}
	\end{center}
\end{table}

We consider selected problems from the NLEVP collection~\cite{betcke2013}, using implementations provided by the NEP-PACK~\cite{jarlebring2018} Julia package. \autoref{tab:parameters} provides a list of our test cases with their parameters.
A Julia implementation of the proposed algorithm is used to solve the problems, giving examples of numerical performance on a variety of NEPs.

A circular contour with trapezoidal quadrature is used for all methods.
In the convergence plots shown, the first iteration corresponds with Beyn's method, allowing us to compare it to the proposed iterative method.
The SVD method is used to solve the linearization, as opposed to the QR approach outlined in algorithm~\ref{alg:hybrid}, to make a fairer comparison with Beyn, which typically uses an SVD.
Eigenvectors are normalized such that \(\norm{x}_2 = 1\) and residuals are then computed as \(\norm{T(\lambda)x}_2/\norm{T(\lambda)}_F\).

A critical aspect of performance in contour based methods is the number of linear system solves.
At each iteration there are \( N \) block linear system solves, where \( N \) is the number of quadrature nodes.
In the problems considered, as long as enough nodes are used to ensure relatively quick convergence, less linear system solves are required in total than for Beyn's method with enough nodes to give comparable convergence.
This is even more significant when the problem is small enough for factorizations of the linear systems to be stored, giving successive iterations a large jump in performance.

\begin{figure}[htbp]
	\centering
	\begin{tikzpicture}
	\begin{axis}[
		ymode=log,
		xmin=0.5,
		xmax=13,
		xtick=data,
		unbounded coords=jump,
		legend style={font=\scriptsize, fill opacity=0.9},
		xlabel={Iteration},
		ymin=5.0e-17,
		ylabel={Residual},
		ymin=5.0e-17,
		extra x ticks={1},
		extra x tick labels={{\scriptsize(Beyn)}},
		extra x tick style={{yshift=-2ex, grid=none}},
		legend entries={\(N=4\), \(N=8\), \(N=16\), \(N=32\), \(N=64\), \(N=128\), \(N=256\) } ]
		\addplot table [x=iter, y=$2$] {bf.dat};
		\addplot table [x=iter, y=$3$] {bf.dat};
		\addplot table [x=iter, y=$4$] {bf.dat};
		\addplot table [x=iter, y=$5$] {bf.dat};
		\addplot table [x=iter, y=$6$] {bf.dat};
		\addplot table [x=iter, y=$7$] {bf.dat};
		\addplot table [x=iter, y=$8$] {bf.dat};
		\fill [blue, fill opacity=0.1] (0.5, 10.0) rectangle (1.5, 5.0e-17);
	\end{axis}
	\end{tikzpicture}
	\begin{tikzpicture}
	\begin{axis}[
		xlabel={Real},
		ylabel={Imaginary},
		xmin=0,
		xmax=2,
		ymin=0,
		ymax=2]
		\draw [dashed, thick] (axis cs:1.0,1.0) circle [radius=0.5];
		\addplot+ [only marks] table {bf_eig_in.dat};
		\addplot+ [only marks, mark=*, mark options={solid, scale=0.35}] table {bf_eig.dat};
	\end{axis}
	\end{tikzpicture}\caption{Butterfly problem.}\label{fig:butterfly}
\end{figure}
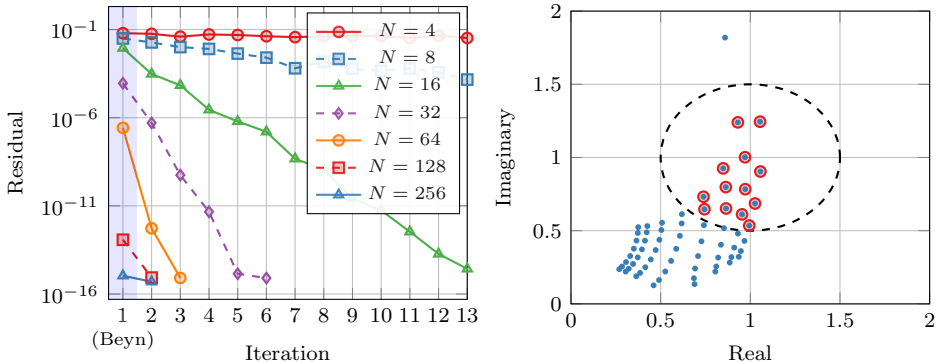

\subsection{Butterfly Problem}
We consider the \verb$butterfly$ problem~\cite{betcke2013}, a quartic eigenvalue problem, as presented by Gavin et.~al.~\cite{gavin2018}.
Choosing the same contour and subspace size \(m\) as used by NLFEAST gives comparable performance. 
As shown in Ref.~\cite{gavin2018a} the subspace size \(m\) can greatly influence convergence speed, in particular with eigenvalues clustered around the contour.
For linear FEAST a rough estimate for a sufficient subspace is \(m = 2 n(\Gamma)\).
Here \(n(\Gamma)=13\) and we intentionally use \(m=20\), slowing convergence in See~\autoref{fig:butterfly}.
Choosing \(m=30\) results in rapid convergence for \(N \ge 16\).

The first iteration corresponds to Beyn's method, where a large number (\(N\ge128\)) of contour nodes are needed for convergence.
Each contour node corresponds to a block linear system solve, the most expensive step of the algorithm.
For this problem the factorizations of the linear systems can easily be stored in memory and reused in subsequent iterations, hence only 16 factorizations are needed with our method.

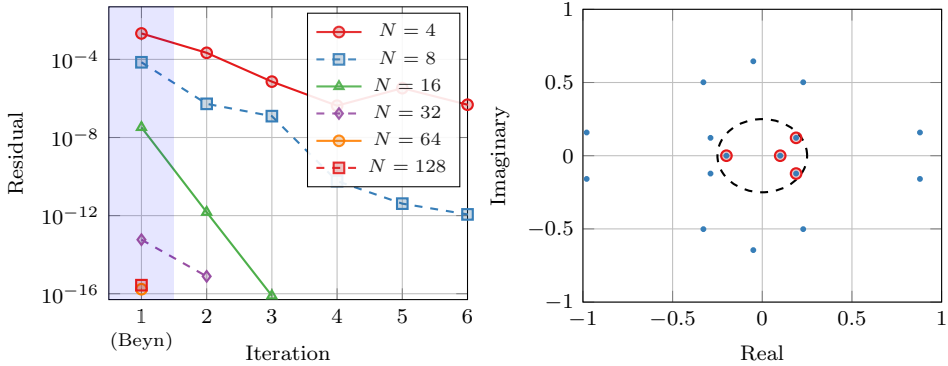
\begin{figure}[htbp]
	\centering
	\begin{tikzpicture}
	\begin{axis}[
		legend style={font=\scriptsize, fill opacity=0.9},
		unbounded coords=jump,
		ymode=log,
		xmin=0.5,
		xmax=6,
		xtick=data,
		xlabel={Iteration},
		ymin=5.0e-17,
		ylabel={Residual},
		extra x ticks={1},
		extra x tick labels={{\scriptsize(Beyn)}},
		extra x tick style={{yshift=-2ex, grid=none}},
		legend entries={\(N=4\), \(N=8\), \(N=16\), \(N=32\), \(N=64\), \(N=128\), \(N=256\) } ]
		\addplot table [x=iter, y=$2$] {deficient.dat};
		\addplot table [x=iter, y=$3$] {deficient.dat};
		\addplot table [x=iter, y=$4$] {deficient.dat};
		\addplot table [x=iter, y=$5$] {deficient.dat};
		\addplot table [x=iter, y=$6$] {deficient.dat};
		\addplot table [x=iter, y=$7$] {deficient.dat};
		\fill [blue, fill opacity=0.1] (0.5, 10.0) rectangle (1.5, 5.0e-17);
	\end{axis}
	\end{tikzpicture}
	\begin{tikzpicture}
	\begin{axis}[
		xlabel={Real},
		ylabel={Imaginary},
		xmin=-1,
		xmax=1,
		ymin=-1,
		ymax=1]
		\draw[dashed, thick] (axis cs:0.0,0.0) circle [radius=0.25];
		\addplot+ [only marks] table {deficient_eig_in.dat};
		\addplot+ [only marks, mark=*, mark options={solid, scale=0.35}] table {deficient_eig.dat};
	\end{axis}
	\end{tikzpicture}\caption{Deficient quadratic problem.}\label{fig:quadratic}
\end{figure}

\subsection{Deficient Quadratic Problem}
This test problem, taken from Beyn~\cite{beyn2012}, illustrates the potential for linearly dependant eigenvectors.
We let \(T_1, T_0 \in \mathbb{C}^{15\times15}\) be random, and set \(T_0\) to have zeroes in the first column, then define the polynomial
\begin{equation}
    T(z) = T_0 + (z-a)(z-b)T_1, \quad a,b \in \mathbb{R}.
\end{equation}
In our tests we set \(a=-0.2, b=0.1\), which are then eigenvalues of \(T\) with the same eigenvector.
It is necessary to use at least two moments, that is to form \(H_0, H_1\) as \(2\times 2\) block moment matrices with algorithm~\ref{alg:generalized}, to resolve these eigenvalues. 
Similarly one can construct a defective matrix polynomial of any degree. 
We use this problem to test our method in the presence of defective eigenvectors.

In~\autoref{fig:quadratic} we see that using \( N=16 \) converges in 3 iterations.
Thus, even without storing factorizations, only 48 linear system solves are needed with our method, where Beyn's method requires \( N=64 \) for similar convergence.

\begin{figure}[htbp]
	\centering
	\begin{tikzpicture}
	\begin{axis}[
		legend style={font=\scriptsize, fill opacity=0.9},
		unbounded coords=jump,
		ymode=log,
		xmax=10,
		xmin=0.5,
		xtick=data,
		xlabel={Iteration},
		ymin=1.0e-18,
		ylabel={Residual},
		extra x ticks={1},
		extra x tick labels={{\scriptsize(Beyn)}},
		extra x tick style={{yshift=-2ex, grid=none}},
		legend entries={\(N=4\), \(N=8\), \(N=16\), \(N=32\), \(N=64\), \(N=128\), \(N=256\) } ]
		\addplot table [x=iter, y=$2$] {gun.dat};
		\addplot table [x=iter, y=$3$] {gun.dat};
		\addplot table [x=iter, y=$4$] {gun.dat};
		\addplot table [x=iter, y=$5$] {gun.dat};
		\addplot table [x=iter, y=$6$] {gun.dat};
		\addplot table [x=iter, y=$7$] {gun.dat};
		\fill [blue, fill opacity=0.1] (0.5, 10.0) rectangle (1.5, 5.0e-19);
	\end{axis}
	\end{tikzpicture}
	\begin{tikzpicture}
	\begin{axis}[
		xlabel={Real},
		ylabel={Imaginary},
		xmin=90000,
		xmax=190000,
		ymin=-50000,
		ymax=50000]
		\draw[dashed, thick] (axis cs:140000.0,0.0) circle [radius=30000.0];
		\addplot+ [only marks] table {gun_eig_in.dat};
		\addplot+ [only marks, mark=*, mark options={solid, scale=0.35}] table {gun_eig.dat};
	\end{axis}
	\end{tikzpicture}
  \caption{Gun problem.}
  \label{fig:gun}
\end{figure}
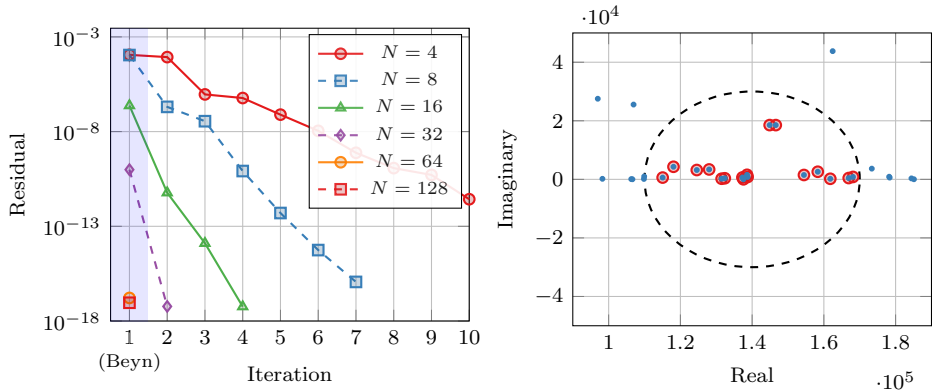

\subsection{Gun Problem}
The \verb$gun$ problem~\cite{betcke2013} models a radio-frequency gun cavity.
The problem is fully nonlinear and has the form
\begin{equation}
    T(\lambda) = K - \lambda M + i \sqrt{\lambda - \sigma_1^2} W_1 + i \sqrt{\lambda - \sigma_2^2} W_2
\end{equation}
where \(M,K,W_1, W_2 \in \mathbb{R}^{9956\times9956}\) are symmetric, and we take \(\sigma_1 = 0, \sigma_2 = 108.8774\).
We use the same contour that Yokota and Sakurai~\cite{yokota2013} choose for demonstrating the SS-RR method.
It is difficult to compare these methods, as SS-RR must solve a projected problem internally with Beyn's method.
Using \(N=32\) the SS-RR method converges to \(10^{-16}\) with \(K=24\) moments and a subspace of size \(m=4\), where with \(K=1\) and \(m=32\) our method converges in two iterations.
See \autoref{fig:gun}.

\begin{figure}[htbp]
	\centering
	\begin{tikzpicture}
	\begin{axis}[
		legend style={font=\scriptsize, fill opacity=0.9},
		unbounded coords=jump,
		ymode=log,
		xmax=10,
		xmin=0.5,
		xtick=data,
		xlabel={Iteration},
		ymin=1.0e-18,
		ylabel={Residual},
		extra x ticks={1},
		extra x tick labels={{\scriptsize(Beyn)}},
		extra x tick style={{yshift=-2ex, grid=none}},
		legend entries={\(N=4\), \(N=8\), \(N=16\), \(N=32\), \(N=64\), \(N=128\), \(N=256\) } ]
		\addplot table [x=iter, y=$2$] {hadeler.dat};
		\addplot table [x=iter, y=$3$] {hadeler.dat};
		\addplot table [x=iter, y=$4$] {hadeler.dat};
		\addplot table [x=iter, y=$5$] {hadeler.dat};
		\addplot table [x=iter, y=$6$] {hadeler.dat};
		\addplot table [x=iter, y=$7$] {hadeler.dat};
		\fill [blue, fill opacity=0.1] (0.5, 10.0) rectangle (1.5, 5.0e-19);
	\end{axis}
	\end{tikzpicture}
	\begin{tikzpicture}
	\begin{axis}[
		xlabel={Real},
		ylabel={Imaginary},
		xmin=-48,
		xmax=-12,
		ymin=-15,
		ymax=15
		]
		\draw[dashed, thick] (axis cs:-30.0,0.0) circle [radius=10.0];
		\addplot+ [only marks] table {hadeler_eig_in.dat};
		\addplot+ [only marks, mark=*, mark options={solid, scale=0.35}] table {hadeler_eig.dat};
	\end{axis}
	\end{tikzpicture}
  \caption{Hadeler problem.}
  \label{fig:hadeler}
\end{figure}
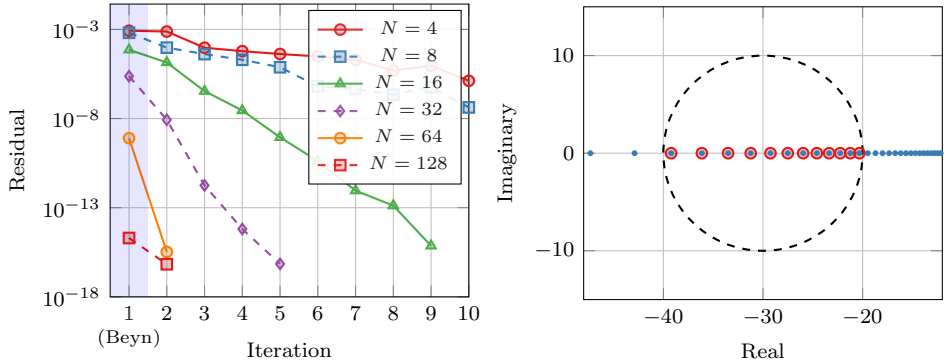

\subsection{Hadeler Problem}
The \verb$hadeler$ problem~\cite{betcke2013} is an NEP of the form
\begin{equation}
    T(\lambda) = (e^\lambda -1)T_2 + \lambda^2 C_1 - \alpha T_0
\end{equation}
where we choose parameters \(\alpha =100\) and \(n=200\). 
The eigenvalues of this problem lie on the real line.
We use the contour chosen in Ref.~\cite{saad2019}.
See~\autoref{fig:hadeler}, we note that even with \(N=128\) Beyn's method is not fully converged.
This data, along with the other problems considered, suggests that using 8 to 32 quadrature nodes is effective in practice, with 16 being a reasonable default choice.

\subsection{Source Code}

The Julia source code for these algorithms is being actively developed and is MIT open source licensed, accessible at \url{https://github.com/spacedome/FEASTSolver.jl}. 


\section{Conclusions}
\label{sec:conclusions}

The method presented in this paper demonstrates numerically the effectiveness of the NLFEAST residual inverse iteration approach applied directly to other contour based algorithms.
Currently this is shown using the two methods presented by Beyn, as an NLFEAST-Beyn hybrid algorithm, further extending this to the SS-type methods would yield a general approach for iterating contour based methods.
A key aspect of this would be extending the RII to incorporate the higher moment expanded subspace, eliminating the need to deflate every iteration.
The flexibility of this iterative method, and the generality of Beyn's Keldysh linearization approach, are steps towards a ``black box'' solver for the nonlinear eigenvalue problem.

The RII has exceptional numerical properties, in linear FEAST it allows for solving the linear systems inexactly, which may extend to the nonlinear RII.
We have found that for some problems it is sufficient to solve the linear systems in single precision, warranting further investigation into the required accuracy of the linear systems solutions for this method, in particular when using iterative solvers, following the work on inexact/iterative linear FEAST~\cite{gavin2017}.

The NLFEAST-Beyn hybrid algorithm gives a practical method for solving NEPs, shown here with numerical validation.
It reduces the number of factorizations needed from Beyn's method, and within a parallel environment, it allows for a significant reduction in the parallel resources used to distribute the contour nodes. In turn, these resources can be reassigned to enable spectrum splitting and/or to solve much larger problems using efficient parallel system solvers.  
Beyn's method is presented with an error analysis based on convergence of the quadrature, but error analysis and theoretical understanding of the NLFEAST-type RII remains an open problem.


\section*{Acknowledgments}
This work is supported by the National Science Foundation, under grants SI2-SSE-1739423 and
 AF-1813480.
\nocite{*}

\bibliographystyle{plain}
\bibliography{references}
\end{document}